\documentclass[12pt]{amsart}
\usepackage{amssymb,amsmath}

\headsep=1cm \numberwithin{equation}{section}
\newtheorem{theorem}{Theorem}[section]
\newtheorem{lemma}[theorem]{Lemma}

\newtheorem{corollary}[theorem]{Corollary}

\theoremstyle{definition}
\newtheorem{definition}[theorem]{Definition}
\theoremstyle{remark}

\newtheorem{example}[theorem]{Example}

\newcommand{\Ass}{\operatorname{Ass}}
\newcommand{\im}{\operatorname{im}}

\newcommand{\pd}{\operatorname{pd}}

\newcommand{\Ext}{\operatorname{Ext}}
\newcommand{\Supp}{\operatorname{Supp}}

\newcommand{\Hom}{\operatorname{Hom}}

\newcommand{\lo}{\longrightarrow}
\newcommand{\fm}{\frak{m}}

\newcommand{\fa}{\frak{a}}

\begin{document}

\begin{center}
{\LARGE\bf On the associated primes of generalized local
cohomology modules }
\\[0.5cm]

{\bf Amir Mafi $^*$}
\\[0.3cm]

{\it Institute of Mathematics, University for Teacher Education,
599 Taleghani Avenue, Tehran 15614, Tehran Iran-and-Department of
Mathematic, Arak University, Arak, Iran}

\end{center}



\thanks{$^*$ Correspondence: Amir Mafi, Institute of Mathematics, University for Teacher Education,
599 Taleghani Avenue, Tehran 15614, Iran; e-mail:
a-mafi@araku.ac.ir.}

\begin{abstract} Let $\fa$ be an ideal of a commutative Noetherian
ring $R$ with identity and let $M$ and $N$ be  two finitely
generated $R$-modules. Let $t$ be a positive integer. It is shown
that $\Ass_R(H_{\fa}^t(M,N))$ is contained in the union of the
sets $\Ass_R(\Ext_R^i(M,H_{\fa}^{t-i}(N)))$, where $0\leq i\leq
t$. As an immediate consequence, it follows that if either
$H_{\fa}^i(N)$ is finitely generated for all $i<t$ or
$\Supp_R(H_{\fa}^i(N))$ is finite for all $i<t$, then
$\Ass_R(H_{\fa}^t(M,N))$ is finite. Also, we prove that if
$d=\pd(M)$ and $n=\dim(N)$ are finite, then $H_{\fa}^{d+n}(M,N)$
is Artinian. In particular, $\Ass_R(H_{\fa}^{d+n}(M,N))$ is a
finite set consisting of maximal ideals.\\
\end{abstract}


\section{Introduction}

Throughout this paper, $R$ is a commutative Noetherian ring with
identity. For an ideal $\fa$ of $R$ and $i\geq 0$, the $i$-th
local cohomology module of $M$ is defined as:
$$H_{\fa}^i(M)=\underset{n}{\varinjlim}\Ext_R^i(R/\fa^n,M).$$
In [{\bf 8}], Huneke conjectured that if $M$ is a finitely
generated $R$-module, then the set of associated primes of
$H_{\fa}^i(M)$ is finite. Singh [{\bf 15}] provides a counter
example for this conjecture. However, it is known that the
conjecture is true in many situations. For example, in [{\bf 11}]
it is shown that if $R$ is local and $\dim R/\fa=1$, then for a
finitely generated $R$-module $M$, the set $\Ass_R(H_{\fa}^i(M))$
is finite for all $i\geq 0$.

Also, Brodmann and Lashgari [{\bf 2}] showed that the first
non-finitely generated local cohomology module of a finitely
generated $R$-module has only finitely many associated primes. Also,
see [{\bf 10}] and [{\bf 4}] for a far reaching generalizations of
this result.

The following generalization of local cohomology theory is due to
Herzog  [{\bf 7}] (see also [{\bf 17}]). The generalized local
cohomology functor $H_{\fa}^i(.,.)$ is defined by
$$H_{\fa}^i(M,N)=\underset{n}{\varinjlim}\Ext_R^i(M/\fa^nM,N)$$
for all $R$-modules $M$ and $N$. Clearly, this is a generalization
of the usual local cohomology functor. Recently, there are some
new interest in generalized local cohomology (see e.g. [{\bf 1}],
[{\bf 5}], [{\bf 6}] and [{\bf 18}]). Our main aim in this paper
is to establish the following.

\begin{theorem} Let $\fa$ be an ideal of $R$ and let $M$ and $N$ be two
finitely generated $R$-modules. Then the following statements hold.
\\(i) For any positive integer $t$,
$$\Ass_R(H_{\fa}^t(M,N))\subseteq
\bigcup_{i=0}^t\Ass_R(\Ext_R^i(M,H_{\fa}^{t-i}(N))).$$\\
(ii) If $d=pd(M)$ and $n=\dim N$ are finite, then
$H_{\fa}^{n+d}(M,N)$ is Artinian. In particular
$\Ass_R(H_{\fa}^{n+d}(M,N))$ consists of finitely many maximal
ideals.\\ (iii) Suppose that $(R,\fm)$ is local with dimension
$n$ and that $d=pd(M)$ is finite.Then
$\Supp_R(H_{\fa}^{n+d-1}(M,N))$ is finite.
\end{theorem}

Clearly (i) extends the main results of [{\bf 10}, Theorem B],
[{\bf 2}, Theorem 2.2] and [{\bf 4}, Corollary 2.7], (ii) extends
[{\bf 12}, Theorem 2.2], and (iii) is an improvement of [{\bf 11},
Corollary 2.4].

\section{The results}

First, we  recall the definition of a weakly Laskerian module. An
$R$-module $M$ is said to be Laskerian if any submodule of $M$ is
an intersection of a finite number of primary submodules.
Obviously, any Noetherian module is Laskerian. In [{\bf 4}], as a
generalization of this notion, we introduced the following
definition.

\begin{definition} An $R$-module $M$ is said to be {\it weakly
Laskerian} if the set of associated primes of any quotient module
of $M$ is finite.
\end{definition}

\begin{example} (i) Every Laskerian module is weakly Laskerian.\\
(ii) Any module with finite support is weakly Laskerian. In
particular, any Artinian $R$-module is weakly Laskerian.
\end{example}

\begin{theorem} Let $\fa$ be an ideal of $R$ and $M$ be  a finitely
generated $R$-module. Let $N$ be an $R$-module and $t$ a positive
integer. Then
$$\Ass_R(H_{\fa}^t(M,N))\subseteq
\bigcup_{i=0}^t\Ass_R(\Ext_R^i(M,H_{\fa}^{t-i}(N))).$$
\end{theorem}

{\bf Proof.} By [{\bf 14}, Theorem 11.38], there is a
Grothendieck spectral sequence
$$E_2^{p,q}:=\Ext_R^p(M,H^q_{\fa}(N))\underset{p}{\Longrightarrow}
H_{\fa}^{p+q}(M,N).$$ For all $i\geq 2$, we consider the exact
sequence
$$0\lo \ker d_i^{0,t}\lo E_i^{0,t}\overset{d_i^{0,t}}\lo
E_i^{i,t-i+1}.              (1)$$\\ Since $E_i^{0,t}=\ker
d_{i-1}^{0,t}/\im d_{i-1}^{1-i,t+i-2}$ and $E_i^{i,j}=0$  for
all  $j< 0$, we may use (1) to obtain  $\ker d_{t+2}^{i,t-i}\cong
E_{t+2}^{i,t-i}\cong \dots =E_{\infty}^{i,t-i}$ for all $0\leq
i\leq t$. There exists a finite  filtration
$$0=\phi^{t+1}H^t\subseteq\phi^tH^t\subseteq\dots\subseteq\phi^1H^t\subseteq\phi^
0H^t=H_{\fa}^t(M,N)$$ such that
$$E_{\infty}^{i,t-i}=\phi^iH^t/\phi^{i+1}H^t$$ for all $0\leq i\leq
t$.

Now, the exact sequences $0\lo \phi^{i+1}H^t\lo \phi^iH^t\lo
 E_{\infty}^{i,t-i}\lo 0$  $(0\leq i\leq t)$
 in conjunction with $$E_{\infty}^{i,t-i}\cong \ker
 d_{t+2}^{i,t-i}\subseteq \ker d_2^{i,t-i}\subseteq E_2^{i,t-i}$$
 yields $$\Ass_R(H_{\fa}^t(M,N))\subseteq
\bigcup_{i=0}^t\Ass_R(\Ext_R^i(M,H_{\fa}^{t-i}(N))). \Box$$

Next, we obtain an extension of [{\bf 2}, Theorem 2.2], [{\bf 10},
Theorem B], and [{\bf 4}, Corollary 2.7].

\begin{corollary} Let $\fa$ be an ideal of $R$, $M$ a finitely generated
$R$-module, and $N$ a weakly Laskerian $R$-module. If
$H_{\fa}^i(N)$ is weakly Laskerian module for all $i<t$, then
$\Ass_R(H_{\fa}^t(M,N))$ is finite.
\end{corollary}

{\bf Proof.} This is immediate by 2.3 and [{\bf 4}, Lemma 2.3 and
Corollary 2.7]. $\Box$

\begin{corollary} Let $(R,\fm)$ be a local ring and let $\dim R\leq 2$.
Let $M$ be a finitely generated $R$- module and $N$ a weakly
Laskerian $R$-module. Then $\Ass_R(H_{\fa}^t(M,N)$ is finite for
all $t\geq 0$.
\end{corollary}

\textbf{Proof.} By [{\bf 11}, Corollaries 2.3, 2.4], [{\bf 3},
Theorem 6.1.2] and 2.2(ii), $H_{\fa}^t(N)$ is weakly Laskerian for
all $t\geq 1$. Also, $\Gamma_{\fa}(N)$  is weakly Laskerian by
2.1. So, by 2.4, $\Ass_R(H_{\fa}^t(M,N))$ is finite for all $t\geq
0$. $\Box$

\begin{corollary} Let $\fa$ be an ideal of a local ring $(R,\fm)$
and $\dim R=n$. Let $M$ be a finitely generated $R$-module and $N$
be an $R$-module such that $H_{\fa}^i(N)=0$ for all $i\neq {n-1},
n$. Then $\Ass_R(H_{\fa}^t(M,N))$ is finite for all $t\geq 0$.
\end{corollary}

\textbf{Proof.} By [{\bf 11}, Corollaries 2.3, 2.4] and
hypothesis, $\Supp_R(H_{\fa}^t(N))$ is finite for all $t\geq 0$;
so that, by 2.4, $\Ass_R(H_{\fa}^t(M,N))$ is finite for all
$t\geq 0$. $\Box$

\begin{corollary} Let $\fa$ be an ideal of a local ring $(R,\fm)$
with $\dim R/\fa=1$. Let $M$ and $N$ be two $R$-modules. Then
$\Ass_R(H_{\fa}^t(M,N))$ is finite for all $t\geq 0$.
\end{corollary}

{\bf Proof.} It is clear that $\Supp_R(H_{\fa}^t(N))$ is finite
for all $t\geq 0$; hence by 2.4, $\Ass_R(H_{\fa}^t(M,N))$ is
finite for all $t\geq 0$. $\Box$\\

 Following [{\bf 9}], a sequence $x_1,\dots ,x_n$
of elements of $\fa$ is said to be an $\fa$-filter regular
sequence on $N$, if
$$\Supp_R((x_1,\dots, x_{i-1})N :_N x_i/(x_1,\dots,
x_{i-1})N)\subseteq V(\fa)$$ for all $i=1,\dots,n$, where $V(\fa)$
denotes the set of all prime ideals of $R$ containing $\fa$. The
concept of an $\fa$-filter regular sequence is a generalization
of the one of a filter regular sequence which has been studied in
[{\bf16}, Appendix 2(ii)] and has led to some interesting results.
 It is easy to see that the analogue of [{\bf16}, Appendix
2(ii)] holds true whenever $R$ is Noetherian, $N$ is a finitely
generated $R$-module and $\fm$ replaced by $\fa$; so that, if
$x_1,\ldots.x_n$ be an $\fa$-filter regular sequence on $N$, then
there is an element $y\in \fa$ such that $x_1,\ldots, x_n,y$ is an
$\fa$-filter regular sequence on $N$. Thus for a positive integer
$n$, there exists an $\fa$-filter regular sequence on $N$ of
length
n.\\

The following Lemma, which needs the concept of a filter regular
sequence, is a generalization of [{\bf 13}, Lemma 3.4].

\begin{lemma} Let $\fa$ be an ideal of $R$ and $M$ be a finitely generated $R$-module
 such that $d=pd(M)$ is finite. Let $N$ be an $R$-module and assume that
$n\in\mathbb{N}$ and $x_1,\dots ,x_n$ is an $\fa$-filter regular
sequence on $N$.Then $ H_{\fa}^{i+n}(M,N)\cong
H_{\fa}^i(M,H^n_{(x_1,\dots ,x_n)}(N))$ for all $i\geq d$.
\end{lemma}

{\bf Proof.} Consider the spectral sequence
$$E_2^{p,q}:=H_{\fa}^p(M,H^q_{(x_1,\dots
,x_n)}(N))\underset{p}{\Longrightarrow} H_{\fa}^{p+q}(M,N).$$ We
have $E_{2}^{p,q}=0$ for $q>n$ (by Theorem 3.3.1 of [{\bf 3}])
and for $q=n$, $p>d$ (by Proposition 2.5 of [{\bf 13}] and Lemma
1.1 of [{\bf 18}]). It therefore follows $E_{2}^{i,n}\cong
E_{\infty}^{i,n}$ and $E_{\infty}^{i,n}\cong H_{\fa}^{i+n}(M,N)$.
This proves the result. $\Box$ \\

The following result is a generalization of [{\bf 12}, Theorem 2.2]
and [{\bf 6}, Theorem 1.2].

\begin{theorem} Let $\fa$ be an ideal of $R$ and let  $M$ and $N$ be two finitely generated
$R$-modules. Assume that $d=\pd(M)$ and $n=\dim N$ are finite.
Then $H_{\fa}^{n+d}(M,N)$ is an Artinian $R$-module. In
particular, $\Ass_R(H_{\fa}^{n+d}(M,N)$ is a finite set consisting
of maximal ideals.
\end{theorem}

{\bf Proof.} Let $x_1,\dots ,x_n$ be an $\fa$-filter regular
sequence on $N$. Then, by 2.8,
$$H_{\fa}^{n+d}(M,N)\cong H_{\fa}^d(M,H^n_{(x_1,\dots ,x_n)}(N))$$
and, by [{\bf 3}, Exercise 7.1.7], $H^n_{(x_1,\dots ,x_n)}(N)$ is
Artinian. Put $S=H^n_{(x_1,\dots ,x_n)}(N)$. Then
$H_{\fa}^d(M,S)\cong H^d(\Hom(M,\Gamma_{\fa}(\dot{E})))$ by [{\bf
6}, Lemma 2.1], where $\dot{E}$ is an injective resolution of $S$
such that its terms are all Artinian modules. Therefore
$H_{\fa}^{n+d}(M,N)$ is Artinian and $\Ass_R(H_{\fa}^{n+d}(M,N))$
is
a finite set consisting of maximal ideals. $\Box$ \\

The following theorem is an improvement of [{\bf 11}, Corollary
2.4].

\begin{theorem} Let $(R,\fm)$ be a local ring of dimension $n$, $N$
an $R$-module, $M$ a finitely generated $R$-module, and $d=pd(M)$
is finite. Then $\Supp_R(H_{\fa}^{n+d-1}(M,N))$ is finite.
\end{theorem}

\textbf{Proof.} Consider the Grothendieck spectral sequence
$$E_2^{p,q}:=\Ext_R^p(M,H^q_{\fa}(N))\underset{p}{\Longrightarrow}
H_{\fa}^{p+q}(M,N).$$ So, we have a finite filtration
$$0=\phi^{d+n}H^{d+n-1}\subseteq \phi^{d+n-1}H^{d+n-1}\subseteq \dots
\subseteq \phi^1H^{d+n-1}\subseteq
\phi^0H^{d+n-1}=H_{\fa}^{d+n-1}(M,N)$$ and the equalities
$E_{\infty}^{i,d+n-i-1}= \phi^iH^{n+d-1}/\phi^{i+1}H^{n+d-1}$ for
all $o\leq i\leq
{n+d-1}$.\\
Since $E_2^{i,n+d-i-1}=0$ for all $i\neq {d-1},d$ and
$E_{\infty}^{i,n+d-i-1}$ is a subquotient $E_2^{i,n+d-i-1} $, it
follows that
$$\phi^{d+1}H^{n+d-1}=\phi^{d+2}H^{n+d-1}=\ldots
=\phi^{d+n}H^{n+d-1}=0$$ and that
$$\phi^{d-1}H^{n+d-1}=\phi^{d-2}H^{n+d-1}=\ldots=\phi^0H^{n+d-1}=H_{\fa}^{n+d-1}(M,N).$$
Now, using the above consequences in conjunction with  [{\bf 11},
Corollaries 2.3, 2.4], it is easy to see that
 $\Supp_R(E_{\infty}^{d,n-1})$ and $\Supp_R(E_{\infty}^{d-1,n})$
 are
finite sets.\\Next, consider the exact sequence
$$0\longrightarrow E_{\infty}^{d,n-1}\longrightarrow
H_{\fa}^{n+d-1}(M,N)\longrightarrow
E_{\infty}^{d-1,n}\longrightarrow 0,$$ to deduce that
$\Supp_R(H_{\fa}^{n+d-1}(M,N))$ is a finite. $\Box$

\end{document}